\documentclass[10pt]{amsart}
\usepackage{amssymb}
\usepackage{amscd}

\def\today{\number\day\space\ifcase\month\or   January\or February\or March\or April\or May\or June\or   July\or August\or September\or October\or November\or December\fi   \number\year}

\theoremstyle{definition}
\newtheorem{thm}{Theorem}[section]
\newtheorem{lem}[thm]{Lemma}
\newtheorem{prp}[thm]{Proposition}
\newtheorem{dfn}[thm]{Definition}
\newtheorem{cor}[thm]{Corollary}

\newtheorem{exa}[thm]{Example}

\newcommand{\beq}{\begin{equation}}
\newcommand{\eeq}{\end{equation}}
\newcommand{\beqr}{\begin{eqnarray*}}
\newcommand{\eeqr}{\end{eqnarray*}}
\newcommand{\bal}{\begin{align*}}
\newcommand{\eal}{\end{align*}}
\newcommand{\bei}{\begin{itemize}}
\newcommand{\eei}{\end{itemize}}

\newcommand{\af}{\alpha}
\newcommand{\bt}{\beta}
\newcommand{\gm}{\gamma}
\newcommand{\dt}{\delta}

\newcommand{\et}{\eta}

\newcommand{\io}{\iota}
\newcommand{\te}{\theta}
\newcommand{\ld}{\lambda}

\newcommand{\ph}{\varphi}
\newcommand{\ps}{\psi}
\newcommand{\rh}{\rho}
\newcommand{\om}{\omega}
\newcommand{\ta}{\tau}

\newcommand{\Ld}{\Lambda}

\newcommand{\Q}{{\mathbb{Q}}}
\newcommand{\Z}{{\mathbb{Z}}}
\newcommand{\R}{{\mathbb{R}}}

\newcommand{\diag}{{\mathrm{diag}}}

\newcommand{\rank}{{\mathrm{rank}}}

\newcommand{\spn}{{\mathrm{span}}}
\newcommand{\card}{{\mathrm{card}}}

\newcommand{\Ker}{{\mathrm{Ker}}}

\newcommand{\andeqn}{\,\,\,\,\,\, {\mbox{and}} \,\,\,\,\,\,}
\newcommand{\QED}{\rule{0.4em}{2ex}}

\newcommand{\ca}{C*-algebra}

\newcommand{\tfae}{the following are equivalent}

\newcommand{\cfn}{continuous function}

\title[Noncommutative torus as transformation group
  C*-algebra]{Realization of a simple higher dimensional noncommutative
  torus as a transformation group C*-algebra}

\author{Benjam\'{\i}n Itz\'{a}-Ortiz}
\author{N.~Christopher Phillips}

\date{21~November 2005}

\address{Centro de Investigaci\'{o}n en Matem\'{a}ticas,
 Universidad Aut\'{o}noma del Estado de Hidalgo,
 Pachuca de Soto, Hidalgo, 42090, M\'{e}xico.}
\email[]{itza@uaeh.edu.mx}
\address{Department of Mathematics, University  of Oregon,
      Eugene OR 97403-1222, USA.}
\email[]{ncp@darkwing.uoregon.edu}

\thanks{Research of the second author
 partially supported by NSF grant DMS 0302401.}

\begin{document}

\begin{abstract}
Let $\te$ be a nondegenerate skew symmetric real $d \times d$ matrix,
and let $A_{\te}$
be the corresponding simple higher dimensional noncommutative torus.
Suppose that $d$ is odd,
or that $d \geq 4$
and the entries of $\te$ are not contained
in a quadratic extension of $\Q.$
Then $A_{\te}$ is isomorphic to the transformation group C*-algebra
obtained from a minimal homeomorphism of a compact connected
one dimensional space locally homeomorphic to the product of
the interval and the Cantor set.
The proof uses classification theory of C*-algebras.
\end{abstract}

\maketitle

\section{Introduction}\label{Sec:Intro}

\indent
Let $\te$ be a skew symmetric real $d \times d$ matrix.
Recall that
the noncommutative torus $A_{\te}$ is by definition~\cite{Rf2}
the universal \ca\  generated by unitaries $u_1, u_2, \dots, u_d$
subject to the relations
\[
u_k u_j = \exp (2 \pi i \te_{j, k} ) u_j u_k
\]
for $1 \leq j, k \leq d.$
(Of course, if all $\te_{j, k}$ are integers, it is not really
noncommutative.
Also,
some authors use $\te_{k, j}$ in the commutation relation instead.
See for example~\cite{Ks2}.)
The algebras $A_{\te}$ are natural generalizations
of the rotation algebras to more generators.
They, and their standard smooth subalgebras,
have received considerable attention.
As just a few examples,
we mention~\cite{Rf1}, \cite{El0}, \cite{Bc}, \cite{RS} and~~\cite{EL}.
In~\cite{Ph2}
(also see the unpublished preprint~\cite{Ph11}),
is it proved that every
simple higher dimensional noncommutative torus is an AT~algebra.

In this paper, we prove that almost every
simple higher dimensional ($d \geq 3$) noncommutative torus
can be realized as the transformation group C*-algebra
obtained from a minimal homeomorphism of a compact connected
one dimensional space.
The minimal homeomorphism is an irrational time map
of the suspension flow of the restriction to its minimal set
of a suitable Denjoy homeomorphism of the circle.
The only exceptional cases are when $d$ is even
and there is a quadratic extension of~$\Q$
which contains all the entries of $\te.$
The proof consists of constructing a homeomorphism,
of the type described,
whose transformation group \ca\  has the same Elliott invariant
as $A_{\te},$
and using the classification results of~\cite{Ln}, \cite{LhP},
and~\cite{Ph2}
(also see the unpublished preprint~\cite{Ph11}).

In the first section,
we prove the result under the assumption that the image of
$K_0 (A_{\te})$ under the unique tracial state of $A_{\te}$
has rank at least~$3.$
In Section~\ref{Sec:Rank2},
we prove that the rank can be~$2$ only when $d$ is even
and there is a quadratic extension of~$\Q$
which contains all the entries of $\te.$
In Section~\ref{Sec:3D},
we give a kind of converse result for the
three dimensional case.

The impetus for this paper came from questions asked
during the second author's talk at the Canadian Operator Symposium
in Ottawa in June 2005,
and both authors are grateful for the invitation
to participate in this conference.
The second author would also like to thank Arkady Berenstein
for discussions which led to the proof of Proposition~\ref{P:Quadratic}.

\section{Construction of the homeomorphisms}\label{Sec:Main}

\indent
Denjoy homeomorphisms of the circle are described
in Section~3 of~\cite{PSS}.
In particular,
if $h_0 \colon S^1 \to S^1$ is a
Denjoy homeomorphisms of the circle with rotation number
$\af \in \R \setminus \Q,$
there is an associated set $Q (h_0) \subset S^1$
of ``accessible points'', defined up to a rigid rotation of the circle,
as in Definition 3.5 of~\cite{PSS}.
The homeomorphism $h_0$ has a unique minimal set $X_0,$
which is homeomorphic to the Cantor set,
and the set $Q (h_0)$ can be thought of as the set of
points at which $S^1$ is ``cut'' to build $h_0 |_{X_0}$
from the rotation $R_{\af}$ by $\af.$
In particular, $Q (h_0)$ consists of the points in $S^1$
which lie on a number $n$ of orbits of $R_{\af},$
with $1 \leq n \leq \infty.$

\begin{dfn}\label{D:RDH}
A {\emph{restricted Denjoy homeomorphism}}
is the restriction of a Denjoy homeomorphism $h_0$ to its
unique minimal set.
The restricted Denjoy homeomorphism is said to
have {\emph{cut number}} $n \in \{ 1, 2, \ldots, \infty\}$
if $Q (h_0)$ consists of exactly $n$ orbits of the
associated rotation on $S^1.$
\end{dfn}

The cut number $n$ is called $n (h_0)$ in~\cite{PSS}.
It depends only on the restriction of $h_0$
to its minimal set $X_0,$
because Theorem~5.3 of~\cite{PSS}
implies that $K_0 (C^* (\Z, X_0, h_0)) \cong \Z^{n + 1}.$

We will make systematic use of the suspension flow of a homeomorphism.
See the introduction to~\cite{It1};
also see II.5.5 and II.5.6 of~\cite{dV}.
We reproduce the definition here.

\begin{dfn}\label{D:SuspF}
Let $g \colon X_0 \to X_0$ be a homeomorphism of
a compact Hausdorff space.
Define commuting actions of $\R$ and $\Z$ on $X_0 \times \R$
by
\[
t \cdot (x, s) = (x, \, s + t)
\andeqn
n \cdot (x, s) = (g^n (x), \, s - n)
\]
for $x \in X_0,$ $s, t \in \R,$ and $n \in \Z.$
Let $X = (X_0 \times \R) / \Z,$
and for $x \in X_0$ and $s \in \R$
let $[x, s]$ denote the image of $(x, s)$ in $X.$
The action of $\R$ on $X_0 \times \R$
descends to an action of $\R$ on $X,$
given by the homeomorphisms $h_t ([x, s]) = [x, \, s + t]$
for $x \in X_0$ and $s, t \in \R,$
called the {\emph{suspension flow}} of $g.$
We refer to $h_t$ as the {\emph{time $t$ map}}
of the suspension flow.
\end{dfn}

We will need the following properties of extensions of dynamical systems.
They are surely well known.
However, we know of no reference for Part~(\ref{POE2})
except for Theorem~A.10 of~\cite{It0}
(although the reverse result, going from $Y$ to $X,$
is Corollary IV.1.9 of~\cite{dV}).
Part~(\ref{POE1}) is in Theorem~A.10 of~\cite{It0}
and also in VI.5.21 of~\cite{dV},
but we give the short proof here for completeness.
The proof of Part~(\ref{POE2})
follows the proof of Theorem~2.6 of~\cite{Wl}.

\begin{lem}\label{P:PropOfExt}
Let $g \colon X \to X$ and $h \colon Y \to Y$
be homeomorphisms of compact Hausdorff spaces,
and let $p \colon X \to Y$
be a surjective map such that $h \circ p = p \circ g.$
(Thus, $g$ is an extension of $h.$)
Let $N \subset Y$ be
\[
N = \{ y \in Y \colon \card ( p^{-1} (y) ) > 1 \}.
\]
Then:
\begin{enumerate}
\item\label{POE1}
If $h$ is minimal and $X \setminus p^{-1} (N)$ is dense in $X,$
then $g$ is minimal.
\item\label{POE2}
If $h$ has a unique ergodic measure $\nu,$
and $\nu (N) = 0,$
then $g$ is uniquely ergodic.
\end{enumerate}
\end{lem}

\begin{proof}
For the first part,
let $K \subset X$ be closed, invariant, and not equal to $X.$
Then $X \setminus K$ is open and nonempty,
so contains a point $x$ of $X \setminus p^{-1} (N).$
By the definition of $N,$ we have $p (x) \not\in p (K).$
Since $p (K)$ is a compact invariant subset of $Y,$
we have $p (K) = \varnothing,$ whence $K = \varnothing.$

Now we prove the second part.
We define a $g$-invariant Borel probability measure $\mu$ on $X$
by $\mu (E) = \nu (p (E \cap [X \setminus p^{-1} (N)])$
for a Borel set $E \subset X.$
Let $\ld$ be any other $g$-invariant Borel probability measure on $X.$
Then $F \mapsto \ld (p^{-1} (F))$
is an $h$-invariant Borel probability measure on $Y,$
whence $\ld (p^{-1} (F)) = \nu (F)$ for every Borel set $F \subset Y.$
In particular, $\ld (p^{-1} (N)) = 0.$
Considering subsets of $X \setminus p^{-1} (N),$
it now follows easily that $\ld = \mu.$
\end{proof}

The following lemma is contained in Proposition~V.2 of~\cite{It0}.
For the convenience of the reader, we give the proof here.
Also, note the relevance of Corollary~2.8 of~\cite{It2},
although it won't in fact be used in the proof.

\begin{lem}\label{L:MinOfSF}
Let $g$ be a Denjoy homeomorphism of $S^1$
with rotation number $\af \in \R \setminus \Q.$
Let $g_0 \colon X_0 \to X_0$
be the restriction of $g$ to its unique minimal set.
Let $t \in \R,$ and let $h_t \colon X \to X$
be the time $t$ map of the suspension flow of $g_0.$
Then \tfae:
\begin{enumerate}
\item
$1,$ $t \af,$ and $t$ are linearly independent over $\Q.$
\item
$h_t$ is minimal.
\item
$h_t$ is uniquely ergodic.
\end{enumerate}
\end{lem}

\begin{proof}
Let $Q = Q (g_0) \subset S^1$ be
the countable set of Definition 3.5 of~\cite{PSS}.
First,
observe that there is a surjection $p_0 \colon S^1 \to X_0$
such that,
with $R_{\af}$ being the rotation by $\af$ on $S^1,$
we have $R_{\af} \circ p_0 = p_0 \circ g_0.$
That is, $g_0$ is an extension of $R_{\af}$
in the sense in Lemma~\ref{P:PropOfExt}.
Moreover, the points in $S^1$ whose inverse images are not unique are
exactly the elements of $Q.$
Let $Y = (S^1 \times \R) / \Z$
be the space of the suspension flow of $R_{\af},$
and let $k_t \colon Y \to Y$ be the time $t$ map of this flow.
Then $h_t$ is an extension of $k_t.$
Let $p \colon X \to Y$ be the extension map.
The set $N$ of points in $Y$
whose inverse images under $p$ are not unique
is $\{ [y, s] \in Y \colon y \in Q \}.$

Define
$f \colon Y \to (\R / \Z)^2$ by
$f ([y, s]) = (y + s \af + \Z, \, s + \Z).$
Then $f \circ k_t \circ f^{-1}$ is the homeomorphism of $(\R / \Z)^2$
given by $(y_1, y_2) \mapsto (y_1 + (t \af + \Z), \, y_2 + (t + \Z)).$

Suppose that
$1,$ $t \af,$ and $t$ are not linearly independent over $\Q.$
Then $f \circ k_t \circ f^{-1}$ has two disjoint nonempty
closed invariant sets $Z_1$ and $Z_2.$
(In fact there are uncountably many.)
So $(f \circ p)^{-1} (Z_1)$ and $(f \circ p)^{-1} (Z_2)$
are disjoint nonempty closed $h_t$-invariant subsets of $X.$
Thus $h_t$ is not minimal.
Since each of these sets carries an invariant Borel probability
measure,
$h_t$ is not uniquely ergodic either.

Now suppose that
$1,$ $t \af,$ and $t$ are linearly independent over $\Q.$
Then $k_t$ is minimal and uniquely ergodic,
with Lebesgue measure as the unique invariant measure.
The set $p_0^{-1} (Q)$ is countable, so that
$X_0 \setminus p_0^{-1} (Q)$ is dense in $X_0.$
Therefore $X \setminus p^{-1} (N)$ is dense in $X.$
It follows from Lemma~\ref{P:PropOfExt}(\ref{POE1})
that minimality of $k_t$ implies minimality of $h_t.$
Moreover, $N$ has measure zero because $Q$ is countable.
So it follows from Lemma~\ref{P:PropOfExt}(\ref{POE2})
that unique ergodicity of $k_t$ implies unique ergodicity of $h_t.$
\end{proof}

\begin{thm}\label{T:ExH}
Let $G_0$ be a finitely generated free abelian group,
let $\om \colon \Z \oplus G_0 \to \R$
be a homomorphism such that $\om (1, 0) = 1$
and $\om (\Z \oplus G_0)$ has rank at least three.
Then there exists a
restricted Denjoy homeomorphism $h_0 \colon X_0 \to X_0$
with cut number $\rank (G_0) - 1,$
and a number $t > 0,$
such that the time $t$ map
$h \colon X \to X$ of the suspension flow of $h_0$
has the following properties:
\begin{itemize}
\item
$h$ is minimal and uniquely ergodic.
\item
$X$ is connected.
\item
There is an isomorphism
$K_0 (C^* ( \Z, X, h)) \cong \Z \oplus G_0$
which sends $(1, 0)$ to $[1],$ sends
$K_0 (C^* (\Z, X, h))_+$
to $\{ 0 \} \cup \{g \in \Z \oplus G_0 \colon \om (g) > 0 \},$
and identifies $\om$ with the map $K_0 (C^*(\Z, X, h)) \to \R$
induced by the unique tracial state
(coming from the unique ergodic measure on $X$).
\item
There is an isomorphism $K_1 (C^*(\Z, X, h)) \cong \Z \oplus G_0.$
\end{itemize}
\end{thm}

\begin{proof}
Set $G = \Z \oplus G_0.$
We identify $G_0$ with $0 \oplus G_0 \subset G$ in the obvious way.

We first claim that there is a direct sum decomposition
$G_0 = H_0 \oplus H_1 \oplus H_2$ with the following properties:
\begin{itemize}
\item
$H_0 \subset \Ker (\om).$
\item
$\om |_{H_1 \oplus H_2}$ is injective.
\item
$\om (H_1) \subset \Q.$
\item
$H_1$ has rank zero or one.
\item
$\om (H_2) \cap \Q = \{ 0 \}.$
\end{itemize}

To prove this, first observe that $\om (G_0)$
is finitely generated and torsion free,
so that $\om |_{G_0}$ has a right inverse
$f_0 \colon \om (G_0) \to G_0.$
Thus there is a direct sum decomposition
$G_0 = H_0 \oplus (f_0 \circ \om) (G_0)$
with $H_0 = G_0 \cap \Ker (\om).$
If $\om (G_0) \cap \Q = \{ 0 \},$
then take $H_1 = \{ 0 \}$ and $H_2 = (f_0 \circ \om) (G_0).$
Otherwise,
let $\pi \colon \R \to \R / \Q$ be the quotient map.
Then $(\pi \circ \om) (G_0)$
is again finitely generated and torsion free,
so that $\pi |_{\om (G_0)}$ has a right inverse
$f_1 \colon (\pi \circ \om) (G_0) \to \om (G_0).$
Set $H_1 = f_0 ( \Ker (\pi |_{\om (G_0)} ))$
and $H_2 = (f_0 \circ f_1 \circ \pi \circ \om) (G_0),$
giving $(f_0 \circ \om) (G_0) = H_1 \oplus H_2.$
We have $\om (H_2) \cap \Q = \{ 0 \}$ because
$\pi |_{(f_1 \circ \pi \circ \om) (G_0)}$ is injective.
Also, $\om (H_1)$ is a nonzero finitely generated subgroup of $\Q,$
and therefore has rank one.
This proves the claim.

Set $m = \rank (H_1 \oplus H_2)$ and $n = \rank (H_0).$
Write
\[
\om (G_0) = \Z \bt_1 + \Z \bt_2 + \cdots + \Z \bt_m,
\]
with $\bt_1, \bt_2, \ldots, \bt_m$ linearly independent over $\Q.$
If $H_1 \neq \{ 0 \},$
then, using the direct sum decomposition
$\om (G_0) = \om (H_1) \oplus \om (H_2),$
we may assume $\bt_m \in \Q.$
We may also obviously assume $\bt_j > 0$ for $1 \leq j \leq m.$
Since $\rank (\om (G)) \geq 3,$
the numbers $1, \bt_1, \bt_2$ must be linearly independent over $\Q.$

Choose $g_1, \ldots, g_m \in H_1 \oplus H_2$
such that $\om (g_j) = \bt_j$
for $1 \leq j \leq m.$
Then $g_1, \ldots, g_m$ form a basis for $H_1 \oplus H_2.$
Further let $g_{m + 1}, \ldots, g_{m + n}$ form a basis for $H_0.$
Choose an integer $N > 1$ so large that
$N^{-n} \bt_j / \bt_1 < 1$ for $2 \leq j \leq m.$

Set
\[
\gm_1 = 0,
\,\,\,\,\,\,
\gm_2 = \frac{\bt_3}{N^n \bt_1},
\,\,\,\,\,\,
\gm_3 = \frac{\bt_4}{N^n \bt_1},
\,\,\,\,\,\,
\ldots,
\,\,\,\,\,\,
\gm_{m - 1} = \frac{\bt_{m}}{N^n \bt_1},
\]
\[
\gm_m = \frac{1}{N^n},
\,\,\,\,\,\,
\gm_{m + 1} = \frac{1}{N^{n - 1}},
\,\,\,\,\,\,
\ldots,
\,\,\,\,\,\,
\gm_{m + n - 1} = \frac{1}{N}.
\]
By the choice of $N,$
we have $\gm_j \in (0, 1)$ for $2 \leq j \leq m + n - 1.$
Let ${\overline{\gm}}_j$ be the image of $\gm_j$ in $S^1 = \R / \Z.$
Set $\af = \bt_2 / (N^n \bt_1),$
and let ${\overline{\af}}$ be its image in $S^1.$
Define $Q \subset S^1$ by
\[
Q = \{ {\overline{\gm}}_j + l {\overline{\af}}
   \colon {\mbox{$1 \leq j \leq m + n - 1$ and $l \in \Z$}} \}.
\]
Then $Q$ is a countable subset of $S^1$ which is invariant under the
rotation $R_{\af}$ by ${\overline{\af}}.$

We now claim that if $j \neq k$ then
$\gm_j - \gm_k \not\in \Z + \af \Z.$
Set
\[
I = \{ 1, \, m, \, m + 1, \, \ldots, \, m + n - 1 \}.
\]
If $j, k \in I,$ then $\gm_j - \gm_k$ is rational and
$0 < | \gm_j - \gm_k | < 1,$ so $\gm_j - \gm_k \not\in \Z + \af \Z.$
If $j, k \not\in I,$ and $\gm_j - \gm_k \in \Z + \af \Z,$
multiply by $N^n \bt_1.$
We get
$\bt_{j + 1} - \bt_{k + 1} \in N^n \bt_1 \Z + \bt_2 \Z.$
This is a linear dependence of
$\bt_1,$ $\bt_2,$ $\bt_{j + 1},$ and $\bt_{k + 1}$ over $\Q,$
a contradiction because $j + 1, \, k + 1 \geq 3.$
Now suppose that $j \in I$ but $k \not\in I,$ and that $j \neq 1.$
Write $j = m + l$ with $0 \leq l \leq n - 1.$
If $\gm_j - \gm_k \in \Z + \af \Z,$
multiply by $N^n \bt_1,$
getting $N^l \bt_1 - \bt_{k + 1} \in N^n \bt_1 \Z + \bt_2 \Z.$
Since $k + 1 \geq 3,$
the numbers $\bt_1,$ $\bt_2,$ and $\bt_{k + 1}$
are linearly independent over $\Q,$
so this is a contradiction.
If instead $j = 1,$
the same procedure would give
$- \bt_{k + 1} \in N^n \bt_1 \Z + \bt_2 \Z,$
a contradiction for the same reason.
This completes the proof the claim.

By Remark~2 in Section~3 of~\cite{PSS},
there exists a Denjoy homeomorphism
$h_0 \colon S^1 \to S^1$ such that $Q (h_0),$
as in Definition 3.5 of~\cite{PSS},
is equal to $Q.$
Let $X_0$ be its unique minimal set.
Let $A_0 = C^* (\Z, X_0, h_0)$ (called $D_{h_0}$ in~\cite{PSS}).
By Proposition~4.2 of~\cite{PSS},
the algebra $A_0$ has a unique tracial state $\ta_0.$
By Theorem~5.3 and Lemma~6.1 of~\cite{PSS}, there is an isomorphism
$\rh_0 \colon \Z^{m + n} \to K_0 (A_0)$
for which, in terms of the standard
generators $\dt_1, \ldots, \dt_{m + n}$ of $\Z^{m + n},$ one has
$(\ta_0)_* (\rh_0 (\dt_1)) = \af,$
$(\ta_0)_* (\rh_0 (\dt_j)) = \gm_j$ for $2 \leq j \leq m + n - 1,$
and $\rh_0 (\dt_{m + n}) = [1].$

We define a different isomorphism
$\rh_1 \colon \Z^{m + n} \to K_0 (A_0)$
as follows.
We set $\rh_1 (\dt_1) = \rh_0 (\dt_m)$
and $\rh_1 (\dt_j) = \rh_0 (\dt_{j - 1})$
for $2 \leq j \leq m,$ and we further set
\[
\rh_1 (\dt_{m + 1})
= \rh_0 (\dt_{m + 1}) - N \rh_0 (\dt_{m}),
\,\,\,\,\,\,
\rh_1 (\dt_{m + 2})
= \rh_0 (\dt_{m + 2}) - N^{2} \rh_0 (\dt_{m}),
\,\,\,\,\,\,
\ldots,
\]
\[
\rh_1 (\dt_{m + n}) = \rh_0 (\dt_{m + n}) - N^{n} \rh_0 (\dt_{m}).
\]
This gives:
\begin{itemize}
\item
$(\ta_0)_* (\rh_1 (\dt_{1})) = 1 / N^n.$
\item
$(\ta_0)_* (\rh_1 (\dt_j)) = \bt_{j} / (N^n \bt_1)$ for
$2 \leq j \leq m.$
\item
$(\ta_0)_* (\rh_1 (\dt_j)) = 0$
for $m + 1 \leq j \leq m + n.$
\end{itemize}

Now take $t = N^n \bt_1.$
Since $1, \bt_1, \bt_2$ are linearly independent over $\Q,$
it is easy to check that $1, t \af, t$
are linearly independent over $\Q.$
So the time $t$ map $h \colon X \to X$ of the suspension flow of $h_0$
is minimal and uniquely ergodic by Lemma~\ref{L:MinOfSF}.
Also, $X$ is connected by Lemma~1.3 of~\cite{It2}.
Let $\mu$ be the unique $h$-invariant Borel probability measure on $X.$
(It is necessarily obtained following Definition~1.8 of~\cite{It1}
from the unique $h_0$-invariant Borel probability measure
$\mu_0$ on $X_0.$)
Let $\ta$ be the corresponding tracial state on $A= C^* (\Z, X, h).$
By Theorem~1.12 of~\cite{It1}, there is an isomorphism
$\ph \colon \Z \oplus K_0 (A_0) \to K_0 (A)$
such that $\ph (1, 0) = [1]$ and
$\ta_* (\ph (0, \et)) = t \cdot (\ta_0)_* (\et)$
for $\et \in K_0 (A_0).$
We now define $\rh \colon G \to K_0 (A)$ on basis elements by
$\rh (1, 0) = \ph (1, 0)$
and $\rh (g_j) = \ph (0, \, \rh_1 (\dt_j))$
for $1 \leq j \leq m + n.$
This defines an isomorphism such that $\rh (1, 0) = [1]$ and
$\ta_* \circ \rh = \om.$
It follows from Theorem 4.5(1) of~\cite{PhX}
that $\et \in K_0 (A)$ is positive if and only if either
$\et = 0$ or $\ta_* (\et) > 0,$ so $\rh$ is an order isomorphism.

Finally, Theorem~1.12 of~\cite{It1} also implies
$K_1 (A) \cong \Z \oplus K_0 (A_0) \cong G$ as abelian groups.
\end{proof}

\begin{thm}\label{T:Realization}
Let $\te$ be a nondegenerate skew symmetric real $d \times d$ matrix.
Let $A_{\te}$ be the corresponding
(higher dimensional) noncommutative torus,
and let $\ta$ be its unique tracial state.
Suppose that $\rank (\ta_* (K_0 (A_{\te}))) \geq 3.$
Then $A_{\te}$ is isomorphic to the crossed product by a minimal
homeomorphism of a compact connected metric space, obtained as the
irrational time map of the suspension flow of
a restricted Denjoy homeomorphism.
\end{thm}

\begin{proof}
We claim that
that for every skew symmetric real $d \times d$ matrix
(nondegenerate or not),
$\Z [1]$ is a direct summand in $K_0 (A_{\te}).$
We prove this by induction on $d.$
The claim is trivially true for $d = 1.$
Suppose it is known for $d,$
and let $\te$ be a skew symmetric real $(d + 1) \times (d + 1)$ matrix.
Let $\te_0$ be the $d \times d$ upper left corner.
Then there is an automorphism $\af$ of $A_{\te_0},$
determined by the requirement that $\af$ multiply each of the standard
unitary generators of $A_{\te_0}$ by a suitable scalar,
such that $A_{\te} \cong C^* (\Z, A_{\te_0}, \af).$
(See Notation~1.1 of~\cite{Ph2} for the explicit formulas.
Also see the unpublished preprint~\cite{Ph11}.)
In particular, $\af$ is homotopic to the identity.
The Pimsner-Voiculescu exact sequence~\cite{PV} therefore splits
into two short exact sequences.
With $\io \colon A_{\te_0} \to A_{\te}$ being the inclusion map,
one of these is
\[
0 \longrightarrow
K_0 (A_{\te_0}) \stackrel{\io_*}{\longrightarrow}
K_0 (A_{\te}) \longrightarrow
K_1 (A_{\te_0}) \longrightarrow 0.
\]
The sequence splits because $K_1 (A_{\te_0})$ is free.
Thus, $K_0 (A_{\te_0})$ is a summand in $K_0 (A_{\te}),$
and the map carries the summand $\Z [1_{A_{\te_0}}]$ in
$K_0 (A_{\te_0})$ to $\Z [1_{A_{\te}}].$
This proves the claim.

Now let $\te$ be a nondegenerate skew symmetric real
$d \times d$ matrix.
Use the claim to write
$K_0 (A_{\te}) = \Z [1] \oplus G_0$
for some subgroup $G_0 \subset K_0 (A_{\te}),$
necessarily isomorphic to $\Z^{2^{d - 1} - 1}.$
Apply Theorem~\ref{T:ExH} with $\ta_*$ in place of $\om,$
obtaining $h \colon X_0 \to X_0$ and $h \colon X \to X$
as there.
Then $h$ is a minimal homeomorphism,
$X$ is a one dimensional compact connected metric space,
and $C^* (\Z, X, h)$ has the same Elliott invariant as $A_{\te}.$
It follows from Theorem~3.5 of~\cite{Ph2}
(also see the unpublished preprint~\cite{Ph11})
that $A_{\te}$ has tracial rank zero in the sense of~\cite{Ln2}
(is tracially~AF in the sense of~\cite{Ln1}; also see~\cite{Ln}),
and it follows from Theorem~4.6 of~\cite{LhP} that
$C^* (\Z, X, h)$ has tracial rank zero.
It is well known that both algebras are simple, separable, nuclear,
and satisfy the Universal Coefficient Theorem.
Therefore Theorem~5.2 of~\cite{Ln} implies that
$A_{\te} \cong C^* (\Z, X, h).$
\end{proof}

We point out that one can use the same methods to match
the Elliott invariants of other \ca s.
For example, let $\te, \gm \in \R$ be numbers such that
$1, \te, \gm$ are linearly independent over $\Q,$
and let $f \colon S^1 \to \R$ be a \cfn.
Let $\af_{\te, \gm, 1, f}$
be the corresponding noncommutative Furstenberg transformation of
the irrational rotation algebra $A_{\te}$ as in
Definition~1.1 of~\cite{OP2}.
The computation of the Elliott invariant follows from
Lemma~1.7 and Corollary~3.5 of~\cite{OP2},
and the proof of Theorem~\ref{T:Realization}
can be applied to find a restricted Denjoy homeomorphism
and a minimal irrational time map $h \colon X \to X$
of its suspension flow such that
$C^* (\Z, X, h)$ has the same Elliott invariant as
$C^* (\Z, A_{\te}, \af_{\te, \gm, 1, f}).$

\section{The rank of the range of the trace}\label{Sec:Rank2}

\indent
In this section,
we determine when $\rank (\ta_* (K_0 (A_{\te}))) = 2.$
This is possible for a simple higher dimensional noncommutative
torus $A_{\te}.$

\begin{exa}\label{E:QuadEx}
Let $\te_0 \in \R \setminus \Q$ satisfy a nontrivial
quadratic equation with integer coefficients, and let
$\te_1, \ldots, \te_n \in (\Q + \Q \te_0) \setminus \Q.$
Then
the tensor product $A = A_{\te_1} \otimes \cdots \otimes A_{\te_n}$
of irrational rotation algebras
is a simple higher dimensional noncommutative torus
such that $\ta_* (K_0 (A)) \subset \Q + \Q \te_0.$
\end{exa}

It seems not to be possible to obtain $A$ as a crossed product in the
manner of Theorem~\ref{T:Realization}.


We give Elliott's description of $\ta_* (K_0 (A_{\te})).$
First, we need some notation.
We regard the skew symmetric real $d \times d$ matrix $\te$
as a linear map from $\Z^d \wedge \Z^d$ to $\R.$
Following~\cite{El0},
if $\ph \colon \Ld^k \Z^d \to \R$
and $\ps \colon \Ld^l \Z^d \to \R$
are linear, we take, by a slight abuse of notation,
$\ph \wedge \ps \colon \Ld^{k + l} \Z^d \to \R$
to be the functional obtained from the alternating functional
on $(\Z^d)^{k + l}$ defined as the antisymmetrization of
\[
(x_1, x_2, \ldots, x_{k + l})
\mapsto \ph (x_1 \wedge x_2 \wedge \cdots \wedge x_k)
      \ps (x_{k + 1} \wedge x_{k + 2} \wedge \cdots \wedge x_{k + l}).
\]
In a similar way, we take
$\ph \oplus \ps \colon \Ld^k \Z^d \oplus \Ld^l \Z^d \to \R$
to be $(\xi, \et) \mapsto \ph (\xi) + \ps (\et).$

\begin{thm}[Elliott]\label{T:Elliott}
Let $\te$ be a skew symmetric real $d \times d$ matrix.
Let $\ta$ be any tracial state on $A_{\te}.$
Then $\ta_* (K_0 (A_{\te}))$ is the range of the
``exterior exponential'',
given in the notation above by
\[
\exp_{\wedge} (\te)
    = 1 \oplus \te \oplus \tfrac{1}{2} \te \wedge \te
           \oplus \tfrac{1}{6} \te \wedge \te \wedge \te
           \oplus \cdots \colon
     {\textstyle{\Ld}}^{\mathrm{even}} \Z^d \to \R.
\]
\end{thm}

\begin{proof}
See~1.3 and Theorem~3.1 of~\cite{El0}.
\end{proof}

\begin{prp}\label{P:Quadratic}
Let $\te$ be a skew symmetric real $d \times d$ matrix.
Suppose that $A_{\te}$ is simple
and $\rank (\ta_* (K_0 (A_{\te}))) = 2.$
Then $d$ is even, and
there exists $\bt \in \R \setminus \Q$
such that every entry of $\te$
is in $\Q + \Q \bt.$
If $d > 2,$ then $\bt$ satisfies a nontrivial quadratic equation
with rational coefficients.
\end{prp}

\begin{proof}
Without loss of generality, $| \te_{j, k} | < 1$ for all $j, k.$
Since $\rank (\ta_* (K_0 (A_{\te}))) = 2,$
there exists $\bt \in \R \setminus \Q$
such that $\ta_* (K_0 (A_{\te})) \subset \Q + \Q \bt.$
Let $u_1, \ldots, u_d$ be the standard unitary generators of $A_{\te}.$
For $j \neq k,$ the elements $u_j$ and $u_k$ generate a subalgebra
isomorphic to $A_{\te_{j, k}},$ which contains a projection $p$ with
$\ta (p) = | \te_{j, k}|.$
Thus $\te_{j, k} \in \Q + \Q \bt.$

We can now write $\te = C + \bt D$ for
skew symmetric $C, D \in M_d (\Q).$

We claim that
simplicity of $A_{\te}$ implies that $D$ is invertible.
First, simplicity implies that $\te$ is nondegenerate,
that is, there is no
$x \in \Q^d \setminus \{ 0 \}$ such that
$\langle x, \, \te y \rangle \in \Q$ for all $y \in \Q^d.$
This is essentially in~\cite{Sl},
and in the form given it appears as Lemma~1.7 and Theorem~1.9
of~\cite{Ph2}.
(Also see the unpublished preprint~\cite{Ph11}.)
Now suppose $D$ is not invertible.
Then there exists $x \in \Q^d \setminus \{ 0 \}$ such that $D x = 0.$
For every $y \in \Q^d$ we then have
\[
\langle x, \, \te y \rangle
= \langle x, \, C y \rangle + \bt \langle x, \, D y \rangle
= \langle x, \, C y \rangle - \bt \langle D x, \, y \rangle.
\]
The first term is in $\Q$ and the second is zero,
contradicting nondegeneracy of $\te.$
This proves the claim.

Corollary~1 to Theorem~6.3 of~\cite{Jc} now implies that $d$ is even.

Now let $d > 2.$
Then $d \geq 4.$
Regard $D$ as a map $\Z^d \wedge \Z^d \to \Q.$
We claim that, as a map $\Ld^4 \Z^d \to \Q,$
we have $D \wedge D \neq 0.$
It is equivalent to prove this with $\Q^d$ in place of $\Z^d.$
Since $\Q$ is a field,
Theorem~6.3 of~\cite{Jc} allows us to assume that
$D = \diag (S, S, \ldots, S)$
with
$S = \left( \begin{smallmatrix} 0 & 1 \\
                   -1 & 0 \end{smallmatrix} \right).$
(There are no zero diagonal blocks since $D$ is invertible.)
Letting $\dt_1, \dt_2, \ldots, \dt_d$ be the standard basis vectors
for $\Q^d,$
a simple calculation now shows that
$(D \wedge D) (\dt_1 \wedge \dt_2 \wedge \dt_3 \wedge \dt_4)
          = \frac{1}{3}.$
This proves the claim.

It remains to show that $\bt$ satisfies a nontrivial quadratic
equation.
By Theorem~\ref{T:Elliott}, the range of
$\frac{1}{2} \te \wedge \te$ is contained in $\Q + \Q \bt.$
Choose $\xi \in \Ld^4 \Z^d$ such that $(D \wedge D) \xi \neq 0.$
Then
\[
(\te \wedge \te) \xi
= (C \wedge C) \xi + \bt (C \wedge D + D \wedge C) \xi
   + \bt^2 (D \wedge D) \xi.
\]
Except for $\bt^2 (D \wedge D) \xi,$
all terms on both sides of this equation
are known to be in $\Q + \Q \bt.$
Since $(D \wedge D) \xi \in \Q \setminus \{ 0 \},$
it follows that $\bt^2 \in \Q + \Q \bt.$
This completes the proof.
\end{proof}

\begin{cor}\label{C:Odd}
Let $\te$ be a nondegenerate skew symmetric real $d \times d$ matrix,
with $d$ odd.
Then $A_{\te}$ is isomorphic to the crossed product by a minimal
homeomorphism of a compact connected metric space, obtained as the
irrational time map of the suspension flow of
a restricted Denjoy homeomorphism.
\end{cor}

\begin{proof}
Combine Theorem~\ref{T:Realization} and Proposition~\ref{P:Quadratic}.
\end{proof}

\begin{cor}\label{C:Even}
Let $\te$ be a nondegenerate skew symmetric real $d \times d$ matrix,
with $d \geq 4$ even.
Suppose the field generated by the entries of $\te$
does not have degree~$2$ over $\Q.$
Then $A_{\te}$ is isomorphic to the crossed product by a minimal
homeomorphism of a compact connected metric space, obtained as the
irrational time map of the suspension flow of
a restricted Denjoy homeomorphism.
\end{cor}

\begin{proof}
Again,
combine Theorem~\ref{T:Realization} and Proposition~\ref{P:Quadratic}.
\end{proof}

\section{The three dimensional case}\label{Sec:3D}

\indent
By Corollary~\ref{C:Odd},
every the odd dimensional noncommutative torus
is isomorphic to the crossed product by a minimal
irrational time map of the suspension flow of
a restricted Denjoy homeomorphism.
For the three dimensional case, there is also a reverse result.

\begin{lem}\label{L:3Ndg}
Let $\te$ be a skew symmetric real $3 \times 3$ matrix,
\[
\te = \left(
\begin{array}{ccc}
     0           &  \te_{1, 2} & \te_{1, 3} \\
     -\te_{1, 2} &  0          & \te_{2, 3} \\
     -\te_{1, 3} & -\te_{2, 3} & 0
\end{array}
\right).
\]
Then $\te$ is nondegenerate
(in the sense used in the proof of Proposition~\ref{P:Quadratic})
if and only if
$\dim_{\Q} (\spn_{\Q} (1, \, \te_{1, 2}, \, \te_{1, 3}, \, \te_{2, 3}))
       \geq 3.$
\end{lem}

\begin{proof}
If
$\dim_{\Q} (\spn_{\Q} (1, \, \te_{1, 2}, \, \te_{1, 3}, \, \te_{2, 3}))
    \leq 2,$
then $\te$ is degenerate by Proposition~\ref{P:Quadratic}.

Now suppose that
$\dim_{\Q} (\spn_{\Q} (1, \, \te_{1, 2}, \, \te_{1, 3}, \, \te_{2, 3}))
       \geq 3.$
Then at least two of $\te_{1, 2},$ $\te_{1, 3},$ $\te_{2, 3}$
are rationally independent.
Suppose $\te_{1, 2}$ and $\te_{1, 3}$ are rationally independent;
the other cases are treated similarly.
Let $x \in \Q^3$ satisfy $\langle x, \te y \rangle \in \Q$
for all $y \in \Q^3.$
We use the formula
\[
\langle x, \te y \rangle
= \te_{1, 2} (x_1 y_2 - x_2 y_1) + \te_{1, 3} (x_1 y_3 - x_3 y_1)
      + \te_{2, 3} (x_2 y_3 - x_3 y_2).
\]
Taking $y = (1, 0, 0),$
we get $- \te_{1, 2} x_2 - \te_{1, 3} x_3 \in \Q,$
whence $x_2 = x_3 = 0.$
Taking $y = (0, 1, 0),$ we then get $\te_{1, 2} x_1 \in \Q,$ whence
$x_1 = 0.$
Thus $x = 0,$ and we have proved that $\te$ is nondegenerate.
\end{proof}

\begin{prp}\label{P:3DDenjoy}
Let $h_0 \colon X_0 \to X_0$ be a restricted Denjoy homeomorphism
with cut number~$2$ (Definition~\ref{D:RDH}),
let $t \in \R,$
and let $h \colon X \to X$ be the time $t$ map of the suspension
flow of $h_0.$
Suppose that $h$ is minimal.
Then $C^* (\Z, X, h)$ is isomorphic to a simple three dimensional
noncommutative torus.
\end{prp}

\begin{proof}
Let $\af \in \R \setminus \Q$ be the rotation number
of a Denjoy homeomorphism of the circle whose restriction to its
minimal set is $h_0.$
Lemma~\ref{L:MinOfSF}
implies that $1,$ $t,$ and $t \af$ are linearly independent over $\Q,$
and that $h$ is uniquely ergodic.
Let ${\overline{\af}}$ be the image of $\af$ in $\R / \Z.$
After a suitable rotation, we may
write the set $Q (h_0) \subset S^1$ of Definition 3.5 of~\cite{PSS} as
\[
Q (h_0)
 = \{ l {\overline{\af}} \colon l \in \Z \}
    \cup \{ {\overline{\gm}} + l {\overline{\af}} \colon l \in \Z \}
\]
for some ${\overline{\gm}} \in \R / \Z,$
and choose $\gm \in \R$ whose image in $\R / \Z$ is ${\overline{\gm}}.$
Let $\ta$ be the unique tracial state on $C^* (\Z, X, h).$
Combining Theorem~5.3 and Lemma~6.1 of~\cite{PSS}
with Theorem~1.12 of~\cite{It1},
we get $K_1 (C^* (\Z, X, h)) \cong \Z^4,$
and we can find a basis for $K_0 (C^* (\Z, X, h))$
consisting of $[1]$ and of three elements $g_1,$ $g_2,$ and $g_3$
such that
$\ta_* (g_1) = t,$ $\ta_* (g_2) = t \af,$ and $\ta_* (g_3) = t \gm.$
Set
\[
\te = \left(
\begin{array}{ccc}
     0           &  t          & t \af      \\
     -t          &  0          & t \gm      \\
     -t \af      & -t \gm      & 0
\end{array}
\right).
\]
Since $1,$ $t,$ and $t \af$ are linearly independent over $\Q,$
Lemma~\ref{L:3Ndg} implies that $\te$ is nondegenerate.
Theorem~\ref{T:Elliott}, together with the fact that
$A_{\te}$ is a simple AT~algebra (Theorem~3.8 of~\cite{Ph2};
also see the unpublished preprint~\cite{Ph11})
implies that
$A_{\te}$ has the same Elliott invariant as $C^* (\Z, X, h).$
Therefore $C^* (\Z, X, h) \cong A_{\te}$ as in the
proof of Theorem~\ref{T:Realization}.
\end{proof}

Generally, however, the \ca s of minimal time $t$ maps of
suspension flows of restricted Denjoy homeomorphisms
are not isomorphic to any noncommutative torus.

\begin{prp}\label{P:WrongDim}
Let $h_0$ be a restricted Denjoy homeomorphism
with cut number $n,$
let $t \in \R,$ and let $h$ be the time $t$ map of
the suspension flow of $h_0.$
If $n + 2$ is not a power of~$2,$
then $C^* (\Z, X, h)$ is not isomorphic
to any higher dimensional noncommutative torus.
\end{prp}

\begin{proof}
We have $K_0 (C^* (\Z, X, h)) \cong \Z^{n + 2}$ as a group
by Theorem~5.3 of~\cite{PSS} and Theorem~1.12 of~\cite{It1},
while $K_0 (A_{\te}) \cong \Z^{2^{d - 1}}$ for any
skew symmetric real $d \times d$ matrix $\te.$
\end{proof}

\begin{prp}\label{P:WrongTrDeg}
Let $d \geq 4.$
Then there exists a restricted Denjoy homeomorphism $h_0$
with cut number $n = 2^{d - 1} - 2,$
and $t \in \R,$
such that the time $t$ map $h$ of
the suspension flow of $h_0$ is minimal,
such that $K_* (C^* (\Z, X, h))$ is isomorphic to the
K-theory of a $d$-dimensional noncommutative torus as a graded
abelian group,
but such that $C^* (\Z, X, h)$ is not isomorphic
to any higher dimensional noncommutative torus.
\end{prp}

\begin{proof}
Choose $t, \af, \gm_1, \gm_2, \ldots, \gm_n \in \R$
with $\gm_1 = 0$ and such that
$t, \af, \gm_2, \gm_3, \ldots, \gm_n$ are algebraically independent
over $\Q.$
Let ${\overline{\gm}}_j$ be the image of $\gm_j$ in $S^1 = \R / \Z,$
and let ${\overline{\af}}$ be the image of $\af$ in $S^1.$
Define $Q \subset S^1$ by
\[
Q = \{ {\overline{\gm}}_j + l {\overline{\af}}
   \colon {\mbox{$1 \leq j \leq n$ and $l \in \Z$}} \}.
\]
Then $Q$ is a countable subset of $S^1$ which is invariant under the
rotation $R_{\af}$ by ${\overline{\af}},$
and, by algebraic independence,
we have $\gm_j - \gm_k \not\in \Z + \af \Z$ for $j \neq k.$
By Remark~2 in Section~3 of~\cite{PSS},
there exists a Denjoy homeomorphism
$h_0 \colon S^1 \to S^1$ such that $Q (h_0),$
as in Definition 3.5 of~\cite{PSS},
is equal to $Q.$
Also write $h_0 \colon X_0 \to X_0$ for the corresponding
restricted Denjoy homeomorphism.
Let $\ta$ be the unique tracial state on $C^* (\Z, X_0, h_0).$
By Theorem~5.3 of~\cite{PSS},
$\ta_* (K_0 (C^* (\Z, X_0, h_0)))$ contains all the numbers
$\af, \gm_2, \gm_3, \ldots, \gm_n.$

Let $h$ be the time $t$ map of the suspension flow.
Then Theorem~1.12 of~\cite{It1} implies that the
range of any tracial state on $K_0 (C^* (\Z, X, h))$
contains all the numbers
$t, \, t \af, \, t \gm_2, \, t \gm_3, \, \ldots, \, t \gm_n.$
By algebraic independence,
this range generates a subfield of $\R$
with transcendence degree at least $n + 1$ over $\Q.$

If $\te$ is any skew symmetric real $d \times d$ matrix,
then Theorem~\ref{T:Elliott}
implies that the image $\ta_* (K_0 (A_{\te}))$
of the K-theory under the trace is contained in the subfield of $\Q$
generated by the entries of $\te,$
which has transcendence degree at most
$\tfrac{1}{2} d (d - 1)$ over $\Q.$
Since $d \geq 4,$
we have $n + 1 > \tfrac{1}{2} d (d - 1),$
so $C^* (\Z, X, h) \not\cong A_{\te}.$

Isomorphism of the K-groups as abelian groups follows from
Theorem~5.3 of~\cite{PSS} and Theorem~1.12 of~\cite{It1}.
\end{proof}

There are surely examples in which an isomorphism can't be
ruled out by transcendence degree,
but can be ruled out by more careful arithmetic.

\end{document}